\documentclass{aip-cp}

\usepackage[numbers]{natbib}
\usepackage{rotating}
\usepackage{graphicx}
\usepackage{verbatim}

\def\ds{\displaystyle}
\def\eps{\varepsilon}

\def\CC{\mathbb{C}}

\def\O{\mathcal{O}}

\newcommand{\eqref}[1]{{(\ref{#1})}}
\newtheorem{theorem}{Theorem}
\newtheorem{lemma}{Lemma}
\newtheorem{corollary}{Corollary} 
\newtheorem{remark}{Remark}

\newcommand{\Comma}{\: ,}
\newcommand{\period}{\: .}

\begin{document}

\title{Optimally truncated WKB approximation for the highly oscillatory stationary 1D Schrödinger equation}

\author[aff1]{Jannis K\"orner\corref{cor1}}

\author[aff1]{Anton Arnold}
\eaddress{anton.arnold@tuwien.ac.at}

\author[aff2]{Christian Klein}
\eaddress{christian.klein@u-bourgogne.fr}

\author[aff1]{Jens Markus Melenk}
\eaddress{melenk@tuwien.ac.at}

\affil[aff1]{Inst.\ f.\ Analysis u.\ Scientific Computing, Technische Universit\"at Wien,  Wiedner Hauptstr. 8, A-1040 Wien, Austria.}
\affil[aff2]{Institut de Math\'{e}matiques de Bourgogne, Universit\'{e} de Bourgogne-Franche-Comt\'{e}, 9 avenue Alain Savary, France.}
\corresp[cor1]{Corresponding author: jannis.koerner@tuwien.ac.at}

\maketitle

\begin{abstract}
We discuss the numerical solution of initial value problems for $\varepsilon^2\,\varphi''+a(x)\,\varphi=0$ in the highly oscillatory regime, i.e., with $a(x)>0$ and $0<\varepsilon\ll 1$. We analyze and implement an approximate solution based on the well-known WKB-ansatz. The resulting approximation error is of magnitude $\mathcal{O}(\varepsilon^{N})$ where $N$ refers to the truncation order of the underlying asymptotic series. When the optimal truncation order $N_{opt}$ is chosen, the error behaves like $\mathcal{O}(\varepsilon^{-2}\exp(-c\varepsilon^{-1}))$ with some $c>0$.
\end{abstract}

\section{INTRODUCTION}

This paper is concerned with the numerical solution of highly oscillatory ordinary differential equations (ODEs) of the form
\begin{equation} \label{schroedinger_eq}
\displaystyle \varepsilon^2 \varphi''(x) + a(x) \varphi(x)= 0 \,, \quad x \in I:=[\xi,\eta]\,;\qquad
\ds \varphi(\xi)=\varphi_0\in\CC\,,\quad
\ds \varepsilon \varphi'(\xi)=\varphi_1\in\CC \,\Comma
\end{equation}
where $0<\eps \ll 1$ is a small parameter and $a(x)\ge a_0>0$ a sufficiently smooth function.
Such initial value problems (IVPs) occur in many applications, e.g.\ in quantum transport \cite{BAP06, Cla_Naoufel}, cosmology \cite{AHLH20}, and mechanical systems (see references in \cite{lub}).

Since, for $\eps\ll 1$, solutions to \eqref{schroedinger_eq} exhibit highly oscillatory behavior, standard ODE-solvers become inefficient as they would have to resolve each oscillation with a step size typically of the size $h=\O(\eps)$. The work presented in \cite{lub} introduces an $\eps$-uniform scheme that achieves accuracy of $\O(h^2)$ for large step sizes up to $h =\O(\sqrt\eps)$. In \cite{ABN,ADK}, schemes of accuracy $\mathcal{O}(h^{2})$ and $\mathcal{O}(h^{3})$ based on a (w.r.t.\ $\eps$) second order WKB approximation of \eqref{schroedinger_eq} were constructed. These methods share the property of \emph{asymptotical correctness}, meaning that the error decreases with $\eps$ even on a coarse spatial grid, if the phase function can be obtained exactly or with spectral accuracy \citep{AKU}. In contrast to these stepping procedures, the approach presented here directly implements a WKB approximation of the solution to (\ref{schroedinger_eq}), including a discussion on the optimal truncation of the underlying asymptotic series; for a detailed analysis and a more extensive presentation of numerical experiments we refer to \cite{AKKM}.


\section{WKB APPROXIMATION}

The well-known WKB-ansatz \cite{bender99,LL85} for the ODE (\ref{schroedinger_eq}) is given by
\begin{equation}\label{WKB-ansatz}
	\varphi(x)\sim \exp\left(\frac{1}{\varepsilon}\sum_{n=0}^{\infty}\varepsilon^{n}S_{n}(x)\right)\Comma\quad \varepsilon \to 0\Comma	
\end{equation}
where $S_{n}$ are unknown complex-valued functions at this point. It is important to note that the asymptotic series in the exponential is typically divergent (as usual for asymptotic series) and must therefore be truncated in order to obtain an approximate solution.

Upon substituting the ansatz (\ref{WKB-ansatz}) into (\ref{schroedinger_eq}), one obtains
\begin{equation}
	\left(\sum_{n=0}^{\infty}\varepsilon^{n}S_{n}^{\prime}(x)\right)^{2}+\sum_{n=0}^{\infty}\varepsilon^{n+1}S_{n}^{\prime\prime}(x)+a(x)=0.
\end{equation}
Then, by comparing powers of $\varepsilon$, we derive a recurrence relation for the functions $S_{n}'$ as follows:
\begin{eqnarray}\label{Sn_recursion}
	S_{0}^{\prime}=\pm i \sqrt{a}\Comma \quad 
	S_{1}^{\prime}=-\frac{S_{0}^{\prime\prime}}{2S_{0}^{\prime}}=-\frac{a^{\prime}}{4a}\Comma \quad
	S_{n}^{\prime}=-\frac{1}{2S_{0}^{\prime}}\left(\sum_{j=1}^{n-1}S_{j}^{\prime}S_{n-j}^{\prime}+S_{n-1}^{\prime\prime}\right)\Comma\quad n\geq 2\period
\end{eqnarray}
Thus, for the computation of each $S_{n}$ for $n\geq 0$ we get one integration constant. 
Additionally, the presence of two different signs in the equation for $S_{0}'$ implies the existence of two sequences of functions satisfying (\ref{Sn_recursion}). This corresponds to the fact that the ODE in (\ref{schroedinger_eq}) has two fundamental solutions. We will denote the sequence induced by the choice $S_{0}^{\prime}=- i \sqrt{a}$ by $(S_{n}^{-})_{n\in\mathbb{N}_{0}}$, and the sequence following from $S_{0}^{\prime}= i \sqrt{a}$ by $(S_{n}^{+})_{n\in\mathbb{N}_{0}}$. It then holds $(S_{2n}^{+})'=-(S_{2n}^{-})'$ and $(S_{2n+1}^{+})'=(S_{2n+1}^{-})'$ for $n\in\mathbb{N}_{0}$. 
%
%
As both sequences $(S_{n}^{\pm})_{n\in\mathbb{N}_{0}}$ provide approximate solutions to the ODE in (\ref{schroedinger_eq}), the general solution can be expressed as a linear combination:
\begin{equation}\label{general_wkb_solution}
	\varphi\approx\varphi_{N}^{WKB}:=\alpha_{N,\varepsilon}\exp\left(\sum_{n=0}^{N}\varepsilon^{n-1}S_{n}^{-}\right)+\beta_{N,\varepsilon}\exp\left(\sum_{n=0}^{N}\varepsilon^{n-1}S_{n}^{+}\right)\period
\end{equation}
Here, $\alpha_{N,\varepsilon}$ and $\beta_{N,\varepsilon}$ are constant w.r.t.\ $x$ but may depend on $\varepsilon$ and $N$. Note that the integration constants in the computation of $S_{n}^{-}$ and $S_{n}^{+}$ can be ``absorbed'' into $\alpha_{N,\varepsilon}$ and $\beta_{N,\varepsilon}$, respectively. Hence, these integration constants can be set to zero without loss of generality. Consequently, we define
\begin{equation}\label{Sn_definition}
	S_{n}^{\pm}(x):=\int_{\xi}^{x}\left(S_{n}^{\pm}\right)'(\tau)\,\mathrm{d}\tau\period
\end{equation}
%
%
%

%
Finally, to ensure that the WKB approximation (\ref{general_wkb_solution}) satisfies the initial conditions in (\ref{schroedinger_eq}), we set
\begin{eqnarray}
	\alpha_{N,\varepsilon}=\frac{\varphi_{0}\left(\sum_{n=0}^{N}\varepsilon^{n}(S_{n}^{+})^{\prime}(\xi)\right)-\varphi_{1}}{\sum_{n=0}^{N}\varepsilon^{n}\left((S_{n}^{+})^{\prime}(\xi)-(S_{n}^{-})^{\prime}(\xi)\right)}\Comma \quad
	\beta_{N,\varepsilon}=\frac{\varphi_{1}-\varphi_{0}\left(\sum_{n=0}^{N}\varepsilon^{n}(S_{n}^{-})^{\prime}(\xi)\right)}{\sum_{n=0}^{N}\varepsilon^{n}\left((S_{n}^{+})^{\prime}(\xi)-(S_{n}^{-})^{\prime}(\xi)\right)}\period
\end{eqnarray}
%

\section{ERROR ANALYSIS OF THE WKB APPROXIMATION}
The results in this section were all proven in \cite{AKKM}. Hereafter, we will use the notation $S_n$ whenever either $S_{n}^{-}$ or $S_{n}^{+}$ could be inserted. Further, we make the assumption throughout this whole section that $S_{0}'$ (and hence $a(x)$) is analytic on a complex neighbourhood $G\subset\mathbb{C}$ of the interval $I$.
\begin{lemma}\label{corollary_Snk}
	There exist constants $K_{1},K_{2}>0$ depending only on $G$ 
	and $S_{0}'$ such that 
	\begin{eqnarray}\label{dSn}
		\| S_{n}\|_{L^{\infty}(I)}&&\leq (\eta-\xi)\| S_{0}'\|_{L^{\infty}(G)}K_{2}^{n}n^{n}\Comma\quad n\in\mathbb{N}_{0}\Comma\label{Sn_est}\\
		\| S_{n}^{(k)}\|_{L^{\infty}(I)}&&\leq
		\| S_{0}'\|_{L^{\infty}(G)}(k-1)!\,K_{1}^{k-1}K_{2}^{n}n^{n}\Comma\quad n\in\mathbb{N}_{0}\Comma\quad k\in\mathbb{N}\period\label{Snk_est}
	\end{eqnarray}
\end{lemma}
\begin{theorem}\label{theorem_wkb_error}
	Let $\varphi$ be the solution of IVP (\ref{schroedinger_eq}). There exist some constants $\varepsilon_{0}\in(0,1)$ and $C>0$ independent of $N$ such that it holds for $\varepsilon\in(0,\varepsilon_{0}]$:
	\begin{eqnarray}\label{wkb_error}
		\|\varphi-\varphi_{N}^{WKB}\|_{L^{\infty}(I)}&&\leq C\| S_{0}'\|_{L^{\infty}(G)}^{2}\left(|\varphi_{0}|\,\| S_{0}'\|_{L^{\infty}(G)}\sum_{n=0}^{N}\varepsilon^{n}K_{2}^{n}n^{n}+|\varphi_{1}|\right)\nonumber
		\exp\left((\eta-\xi)\| S_{0}'\|_{L^{\infty}(G)}\sum_{n=0}^{\lfloor \frac{N}{2}\rfloor}\varepsilon^{2n}K_{2}^{2n+1}(2n+1)^{2n+1}\right)\\
		&&\quad\cdot\left(2\varepsilon^{N}K_{2}^{N+1}(N+1)^{N+1}+\sum_{n=2}^{N}\varepsilon^{N+n-1}K_{2}^{N+n}\sum_{k=2}^{N+n}k^{k}(N+n-k)^{N+n-k}\right)\Comma
	\end{eqnarray}
	with $K_{2}>0$ being the constant from Lemma \ref{corollary_Snk}.
	In particular, we have $\|\varphi-\varphi_{N}^{WKB}\|_{L^{\infty}(I)}=\mathcal{O}(\varepsilon^{N})$, as $\varepsilon\to 0$.
\end{theorem}
According to Theorem \ref{theorem_wkb_error}, the error of the WKB approximation is $\mathcal{O}(\varepsilon^{N})$ as $\varepsilon \to 0$, for a fixed $N\geq 0$. However, in practical applications the situation is the opposite, namely, the small model parameter $\varepsilon$ is fixed and $N$ can be chosen freely. The following corollary states that the approximation error is exponentially small w.r.t.\ $\varepsilon$, if $N$ is chosen adequately.
\begin{corollary}\label{corollary_optimal}
	Let $\varphi$ be the solution of IVP (\ref{schroedinger_eq}). Then there exists $N=N(\varepsilon)\in \mathbb{N}$ such that
	\begin{equation}\label{optimal_error}
		\| \varphi-\varphi_{N}^{WKB}\|_{L^{\infty}(I)} \leq \frac{C_{1}}{\varepsilon}\exp\left(-\frac{1}{e K_{2}\varepsilon}\right)+\frac{C_{2}}{\varepsilon^{2}}\exp\left(-\frac{\ln(e/2)}{e K_{2}\varepsilon}\right)\Comma
	\end{equation}
	with constants $C_{1},C_{2}>0$ independent of $\varepsilon$. $K_{2}$ is again the constant from Lemma \ref{corollary_Snk}.
\end{corollary}
\begin{remark}\label{remark}
	The proof of Corollary \ref{corollary_optimal} (see \cite{AKKM}) indicates that, in order to obtain estimate (\ref{optimal_error}), $N$ must be chosen proportional to $\varepsilon^{-1}$. Moreover, since the asymptotic series in (\ref{WKB-ansatz}) is typically divergent, there exists an optimal truncation order $N_{opt}$ that minimizes the corresponding approximation error. Then, by definition of $N_{opt}$, the ``optimal error'' $\| \varphi-\varphi_{N_{opt}}^{WKB}\|_{L^{\infty}(I)}$ also satisfies estimate (\ref{optimal_error}) and is hence at most of order $\mathcal{O}(\varepsilon^{-2}\exp(-c/\varepsilon))$ as $\varepsilon\to 0$, where $c>0$. 
\end{remark}
%

\section{NUMERICAL TEST}

\begin{figure}
	\centering
	\includegraphics[scale=0.8]{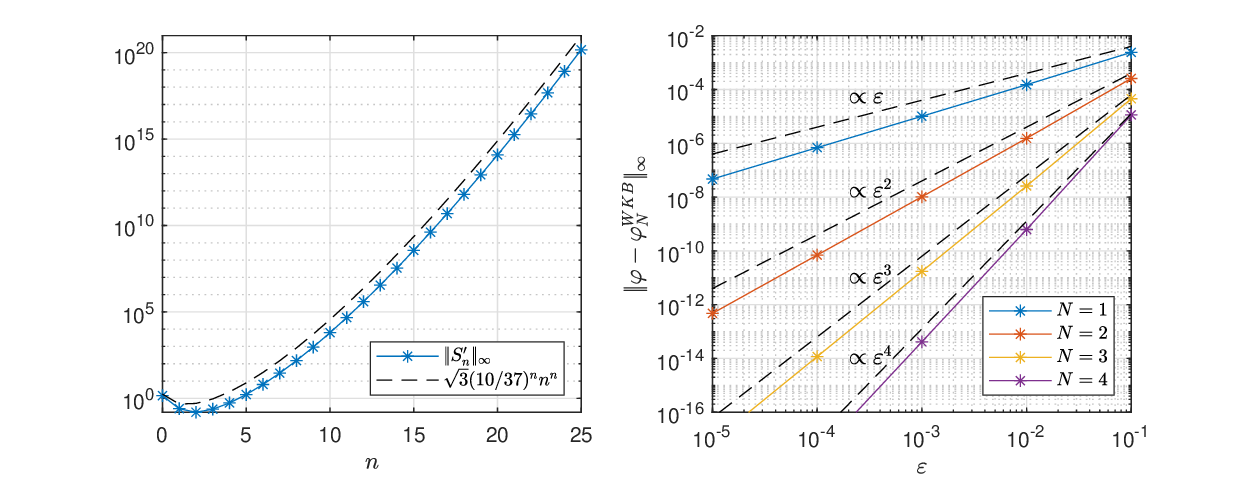}
	\caption{Left: Semilog plot of the $L^{\infty}$-norms of $S_{n}'$ as a function of $n$. Right: Log-Log plot of the $L^{\infty}$-error of the WKB approximation as a function of $\varepsilon$, for several choices of $N$.}
	\label{plot:dSn_norms}
\end{figure}

\begin{figure}
	\centering
	\includegraphics[scale=0.8]{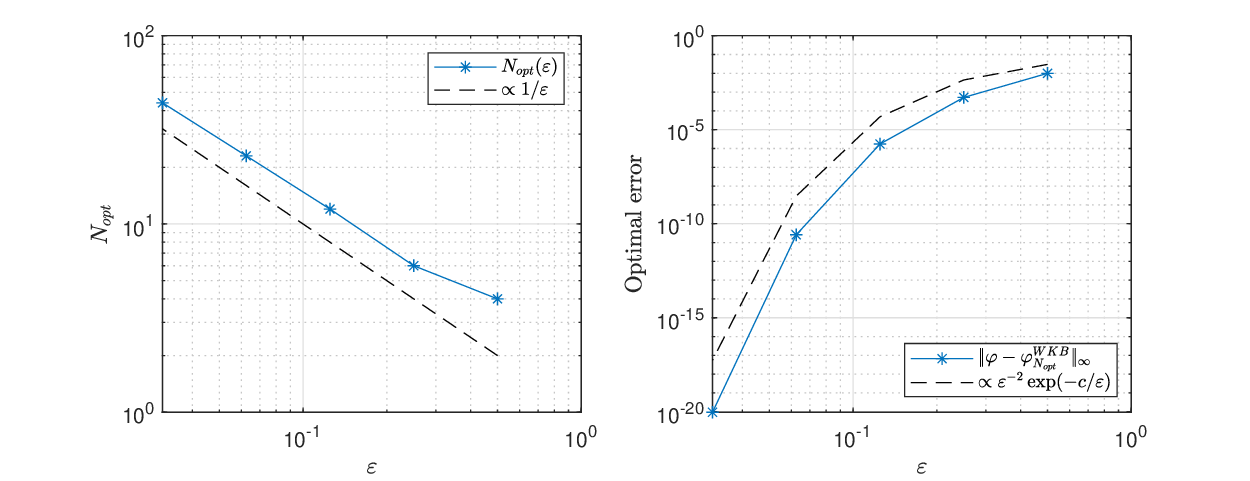}
	\caption{Left: Log-Log plot of the optimal truncation order as a function of $\varepsilon$. Right: The corresponding values of the optimal error as a function of $\varepsilon$, also as a Log-Log plot. The dashed line is proportional to $\varepsilon^{-2}\exp(-\frac{6}{5\varepsilon})$.}
	\label{plot:airy_N_opt_opt_err}
\end{figure}

We investigate the Airy equation on the interval $I=[1,2]$ by choosing $a(x)=x$. The initial conditions for (\ref{schroedinger_eq}) are chosen as $\varphi_0=Ai(-\frac{1}{\varepsilon^{2/3}}) + i Bi(-\frac{1}{\varepsilon^{2/3}})$ and $\varphi_1=-\varepsilon^{1/3}\left(Ai^{\prime}(-\frac{1}{\varepsilon^{2/3}}) + i Bi^{\prime}(-\frac{1}{\varepsilon^{2/3}})\right)$, where $Ai$ and $Bi$ denote the Airy functions of first and second kind, respectively. Since we are dealing with very small approximation errors, we use the Advanpix Multiprecision Computing Toolbox for MATLAB \cite{advanpix} with quadruple-precision to avoid rounding errors. To compute the WKB approximation, we first symbolically calculate the functions $S_n'$ using (\ref{Sn_recursion}). These functions are then integrated using the highly accurate Clenshaw-Curtis quadrature \cite{clencurt60}. In Figure \ref{plot:dSn_norms} we present on the left the $L^\infty$-norm of the functions $S_{n}'$ as a function of $n$. The plot shows good agreement with estimate \eqref{dSn}, as indicated by the dashed line, where we plotted $\sqrt{3}K_{2}^nn^n$ with the experimentally determined constant $K_{2}=10/37$. On the right of Figure \ref{plot:dSn_norms} we plot the $L^{\infty}$-error of the WKB approximation as a function of $\varepsilon$ (with double logarithmic scale), for several choices of $N$. The plot illustrates that the error is indeed of order $\mathcal{O}(\varepsilon^{N})$, as $\varepsilon\to 0$, in agreement with Theorem \ref{theorem_wkb_error}. In Figure \ref{plot:airy_N_opt_opt_err}, on the left, we show the optimal truncation order $N_{opt}$ as a function of $\varepsilon$ and observe that $N_{opt}\propto \varepsilon^{-1}$. On the right, we plot the corresponding optimal error as well as the function $\frac{1}{5\varepsilon^{2}}\exp(-\frac{6}{5\varepsilon})$, corresponding to the second term in (\ref{optimal_error}). We observe good agreement between the plot and the statements of Corollary \ref{corollary_optimal} and Remark \ref{remark}, respectively. 


\section{ACKNOWLEDGMENTS}
The authors JK, AA, and JMM acknowledge support by the project SFB \# F65 of the FWF, and the authors AA and CK also by the bi-national FWF-project I3538-N32.


\nocite{*}
\bibliographystyle{aipnum-cp}%

\end{document}